\documentclass[11pt,eqno,sumlimits]{amsart}

\usepackage{amssymb, amscd, amsmath, epsfig, mathtools}
\usepackage[utf8]{inputenc}
\usepackage{amsthm}
\usepackage{enumerate}
\usepackage{xcolor}
\usepackage{scalerel}
\usepackage{soul}
\usepackage{tikz-cd}
\usepackage{esint}
\usepackage{slashed}
\usepackage[textsize=tiny]{todonotes}
\usepackage{pinlabel} % for labelling the figure
\usepackage{hyperref}
\usepackage{url}
\usepackage{mathtools}
\usepackage{ulem}
\usepackage{ulem}

\usepackage{cancel}

\usepackage[margin=1.0in]{geometry}

\newtheorem{theorem}{Theorem}[section]
\newtheorem*{theorem*}{Theorem}

\newtheorem{corollary}{Corollary}[section]
\newtheorem{lemma}{Lemma}[section]
\newtheorem{proposition}{Proposition}[section]
\theoremstyle{definition}
\newtheorem{remark}{Remark}[section]
\newtheorem{example}{Example}[section]

\newcommand{\C}{\mathbb C}
\newcommand{\R}{\mathbb R}

\newcommand{\Z}{\mathbb Z}
\newcommand{\Q}{\mathbb Q}
\newcommand{\calC}{\mathcal C}

\newcommand{\dvol}{ d\text{Vol}_{g}}

\newcommand{\di}{\slashed{\partial}}

\newcommand{\tg}{\tilde g}
\DeclareMathOperator{\ind}{ind}
\newcommand{\bg}{\bar g}
\newcommand{\sgn}{\text{sgn}}

\newcommand{\xitilde}{\tilde{\xi}}

\DeclareMathOperator{\U}{U}
\newcommand{\spin}{\,\operatorname{spin}}
\newcommand{\spincs}{\mathfrak s}  %% for spin structures
\DeclareMathOperator{\Tr}{Tr}

\newenvironment{introtheorem}[1]{%
  \renewcommand\thetheorem{#1}
  \theorem
}{\endtheorem}

\begin{document}

\title[The conformal Laplacian and P{S}{C} metrics on manifolds with boundary]{The conformal Laplacian and positive scalar curvature metrics on manifolds with boundary}

\author[S. Rosenberg, D. Ruberman and J. Xu]{Steven Rosenberg, Daniel Ruberman and Jie Xu}
\address{
Department of Mathematics and Statistics, Boston University, Boston, MA, USA}
\email{sr@math.bu.edu}
\address{Department of Mathematics, Brandeis University, Waltham, MA, USA}
\email{ruberman@brandeis.edu}
\address{
Institute for Theoretical Sciences, Westlake University, Hangzhou, Zhejiang Province, China;
Department of Mathematics, Northeastern University, Boston, MA, USA}
\email{jie.xu@northeastern.edu}

\thanks{The second author was partially supported by NSF Grant DMS-1928930 while he was in residence at the Simons Laufer Mathematical Sciences Institute (formerly known as MSRI), as well as NSF FRG Grant DMS-1952790.\\
Math.~Subj.~Class.~2020: 53C21 (primary), 35J66, 53C27, 58J05, 58J20 (secondary).}

\date{}							% Activate to display a given date or no date

\begin{abstract}
We give examples of spin $4$-manifolds with boundary $(M,\partial M)$ such that the boundary $\partial M$ has a positive scalar curvature metric which cannot be extended to a positive scalar curvature metric on $M$ with 
mean convex boundary. These manifolds have the equivalent analytic property that for any metric $g$ on $M$, the conformal Laplacian on $M$ with appropriate boundary conditions 
and the conformal Laplacian on $\partial M$ cannot both be positive.
The obstruction to the positivity of the conformal Laplacians is given by a real-valued $\xi$-invariant associated to the APS theorem for the twisted Dirac operator.
We use analytic techniques related to the prescribed scalar curvature problem in conformal geometry to directly treat metrics which are not a product near the boundary.
\end{abstract}

\maketitle

\section{Introduction}

There is a well developed theory of determining which closed manifolds $M$ have topological obstructions to admitting positive scalar curvature (psc) metrics (see \cite{Stolz}).   The corresponding theory for manifolds with boundary 
$(M, \partial M)$ is growing but less complete.   
 The main results 
in this paper are both geometric and analytic extensions of  known results for psc metrics which are a product near the boundary.  On the geometric side, we 
weaken the product metric assumption to the 
mean convex case.  On the analytic side,
we prove results relating psc metrics to the positivity of the conformal Laplacians 
on $M$ and $\partial M$ with
certain geometric  boundary conditions.  The geometric and analytic results are related: 
there exists a metric $g$ on $M$ with positive conformal Laplacians on $M$ and $\partial M$
 iff there exists a psc metric $k$ on $\partial M$ which extends to a psc metric $k'$ on $M$ which is a product near $\partial M$ (Cor.~\ref{cor:1}).
In our approach, the main technical work in the paper is in geometric analysis, specifically in conformal geometry on manifolds with boundary.

To state the main results, let $R_g$ be the scalar curvature of a Riemannian metric $g$ on $(M, \partial M)$ and let $A_g$ be second fundamental form of $\partial M$ with respect to the outward pointing unit normal vector $\nu$. $\partial M$ is totally geodesic if $A_g \equiv 0$,
minimal if 
$h_g :={\rm Tr}(A_g) \equiv 0$, and mean convex if $h_g\nu$ does not point inward anywhere on the boundary 
({\it i.e.}, $h_g\geq 0$).
Let $\eta_1$ be the first eigenvalue of the conformal Laplacian $\Box_g$ with boundary conditions given in (\ref{Robin}), and for the inclusion $\imath:\partial M\to M$, let $ \zeta_{1} $ be the first eigenvalue of the conformal Laplacian $ \Box_{\imath^{*}g} $ for the closed manifold $ (\partial M, \imath^{*}g) $. 
For $M$ a spin $4k$-manifold, let
$\xi_\alpha(\partial M, \imath^*g)$ be the $\R$-valued $\xi$-invariant of $\partial M$ associated to the Dirac operator on $\partial M$, with $\alpha$ a unitary representation of $\pi_1(\partial M).$  

For notation, let $[g]$ denote the conformal class of a metric $g$.
Write $k\in [psc]$ if a metric $k$ is in the conformal class of a psc metric on $\partial M$, and write $k'\in [psc, mc]$ if $k'$ is in the conformal class of a psc metric on $M$ with mean convex boundary.

\begin{introtheorem}{3.1}  Let $(M,\partial M, g)$ be a spin manifold with $\eta_1>0$ and $ \zeta_{1} > 0 $. 
Assume that $\alpha$ is a unitary representation  of $\pi_1(\partial M)$ that extends to a unitary representation of $\pi_1(M).$ Then $\xi_\alpha(\partial M, \imath^*g) = 0.$  
In particular, if $M$ admits a metric $\tg \in [g] $
with $R_{\tg}>0$, mean convex boundary, and $R_{\imath^*\tg}>0$, then  
$\xi_\alpha(\partial M, \imath^*g) = 0.$ 
\end{introtheorem}
  
The $\xi$-invariant is ``quasi-topological," in that its $\R/\Z$ reduction is independent of the metric on $\partial M.$ Using this invariant, we give examples of $4$-manifolds with obstructions to extending psc metrics from their boundary to the interior.
Let $L(p,q)$ be a $3$-dimensional lens space in the usual notation, and let $p^2 L(p,q)$ be $p^2 $ disjoint oriented copies of $L(p,q)$.  

\begin{introtheorem}{4.1}
Let $p$ be odd. Then there exists a $4$-manifold $M(p,q)$ with boundary $p^2L(p,q)$ such that $M(p,q)$ has no metric with $\eta_1>0$ and $\zeta_1>0$. 
 In particular, $p^2L(p,q)$ has no metric $k\in [psc]$ which extends to a metric $k'\in [psc,mc]$ on $M(p,q)$. 
\end{introtheorem}

\noindent   
We believe these are the first examples of these analytic/geometric restrictions in this dimension.  
The construction of $M(p,q)$ is explicit; bordism theory shows the existence of such a manifold with boundary $p L(p,q)$.

Schematically, Thm.~3.1 is of the form  
\begin{align}\label{!} R_g > 0 \ &\text{plus}\ ( \tg\ \text{is  mean convex and}
\ R_{\imath^{*} \tilde{g}} > 0) 
\Rightarrow \eta_1 >0, \zeta_1 >0 \\
&\Rightarrow 
\xi_\alpha(\partial M, \imath^*g)= 0.\notag
\end{align}
The first implication follows 
from the variational definitions of $\eta_1$ and $\zeta_1$,
while the second implication uses 
\cite[Thm.~1.2]{XU2} and the APS index theorem. While the most appealing hypotheses to geometers involve the scalar curvature,   conformal analytic assumptions like $\eta_1>0, \zeta_1>0$ are  easier to work with analytically.

There is a  series of implications similar to (\ref{!}) for 
 closed Riemannian spin manifolds $(M^{4k},g)$.  The Lichnerowicz formula and the AS index theorem give the standard obstruction to psc metrics on $M$:
 $$R_g >0 \Rightarrow 
{\rm ker} (\di^{\pm}_g) = 0 \Rightarrow \hat A(M) = 0,$$ 
 where $\di_g^{\pm}$ is the Dirac operator on the sum of plus and minus spinors.  This  has no immediate connection to conformal geometry.
We can insert $\eta_1$ somewhat artificially: for $\tilde g$ conformally equivalent to $g$ with $R_{\tilde g}$ constant, 
$$ R_g >0 \Rightarrow \eta_1>0\ {\rm and}\ \Delta+ \frac{R_{\tilde g}}{4} >0\Rightarrow
{\rm ker} (\di_{\tilde g}^{\pm}) = {\rm ker} (\di_g^{\pm}) =0 \Rightarrow \hat A(M) = 0.$$
The first implication uses the solution of the Yamabe conjecture: $R_g>0$ implies $\tilde g$ exists
with $R_{\tilde g}$ a positive constant.
For the second implication, see \cite[Thm.~3.5]{R2}. 
In contrast, the use of 
 conformal analytic techniques seems crucial in the manifold with boundary case.

For related work, \cite{BG} uses $\xi$-invariants to produce examples of psc metrics $g_1,g_2$ on boundary components of a spin manifold $M$ such that $M$ admits no psc metric which is a product metric $g_i\oplus dt^2$ near the boundary,     and \cite{XZD} uses $\xi$-invariants as an obstruction to the filling problem in conformally flat ({\it e.g.,} hyperbolic) geometry. 
 \cite{AB2} discusses when a psc metric on  $\partial M$ can be filled in to a psc metric on a manifold $M$ with product metric near $\partial M$, 
 \cite{AB1}  gives conditions under which a psc metric has a deformation to a psc metric which is a product near the boundary (which we use in \S5), and \cite{BH} prove that a psc metric with mean convex boundary can be deformed through psc metrics to a totally geodesic boundary. 
   The product metrics 
in \cite{AB2, AB1, BG}            
have 
mean convex boundary, so the  work in \S3 can be thought of as 
the geometric analysis needed to directly treat  the non-product metric case.
Using a combination of Ricci flow and other dimension-specific techniques, \cite{CL}  completely characterizes 
the topology of all $3$-manifolds with boundary admitting metrics with $R_g>0$ and non-negative mean curvature. 
In a different approach, \cite[Thm.~19]{BH} gives a differential topological obstruction in terms of $K$-area to a manifold with  $R_g>0$ and non-negative mean curvature, 
and produces an example of a $3$-manifold admitting no such metric. 
In the other direction, the existence of constant psc metrics on manifolds with boundary goes back to 
 \cite{E,Gajer} (with \cite{Gajer} corrected in \cite{walsh:psc-I}), with recent work by the third author \cite{ XU3, XU4, XU1}. \cite[Thms.~5.3, 5.5]{RW} gives topological conditions under which a totally non-spin $M$ has a psc metric with vanishing or positive mean curvature.

As an outline of the paper, in \S2  we give the analytic setup.  As a first conformal analytic/geometric result,
we prove that $\eta_1>0, \zeta_1>0$
implies the existence of a metric $\tg$ conformal to $g$ with $R_{\tg}, R_{\imath^*\tg}$ positive (Thm.~\ref{Scalar:thm1}). We also prove that $\eta_1>0$ gives a $\tg$ conformal to $g$ with $R_{\tg}>0,$ 
 $h_{\tg} =0$
(Cor.~\ref{APS:cor1}). Note that the converse of this result is elementary. 
In \S3, we prove Thm.~\ref{thm:4.1}. In \S4, we prove Thm.~\ref{thm:5.1}. 
In \S5, we reprove the main theorems using deep results in the literature.  
Namely, using \cite[Thm.~0.2]{LZ}, 
  we prove that the sign of $\eta_1$ and the associated relative Yamabe invariant of \cite{AB2} are the same (Prop.~\ref{prop:RYI}). Using \cite[Cor.~B]{AB2}, we prove Cor.~\ref{cor:1} mentioned above.  With this equivalence of analytic and geometric quantities, the proofs of the main theorems are quick.

We would like to thank John Lott for his help.

\section{Positivity of the Conformal Laplacian and Positive Scalar Curvature}

Throughout this paper, let $ (M^n, g) $ be a connected, compact, orientable smooth Riemannian manifold with 
smooth boundary $ \partial M $, $ n = \dim M \geqslant  3 $.  It is standard that if $R_g>0$ on $M$ and $h_g = 0$ on $\partial M$, then $\eta_1>0.$  In this section, we consider a series of converse results. For example,  if $\eta_1>0$ and $\zeta_1>0$, there exists a metric $\tg\in [g]$,  such that  the scalar curvatures $R_{\tg}$ on $M$ and
$R_{\imath^*\tg}$ on $\partial M$ are positive constants (Cor~\ref{Ric:cor1}). We also prove the known or folklore result that if $\eta_1>0$, there exists  $ \tilde{g}\in [g]$  such that $ R_{\tilde{g}}>0 $ and $ h_{\tilde{g}} = 0 $  (Cor.~\ref{APS:cor1}).  These results are used in Thm.~\ref{thm:4.1} and Thm.~\ref{thm:5.1}.

\subsection{Notation and basic setup}\label{positive:def1}

The unit outward normal vector field along $ \partial M $ with metric $ g $ is denoted by $\nu = \nu_{g}.$ $ R_{g} $ is the scalar curvature of $ g $,  
and $ A_g $  is the second fundamental form on $ \partial M $: $ A_{g}(X, Y) = \left( \nabla_{X} Y \right)^{\perp} $ for  $ X, Y \in T(\partial M) $, where $ \nabla $ is the Levi-Civita connection 
of $g$ and $\perp$ is projection to the outward normal direction.  The mean curvature  at $ x \in \partial M$ is 
  $ h_{g} = {\rm Tr}(A_g) = g\left(\sum_i \nabla_{e_i} e_i,\nu\right)$, where $\{e_i\}$ is an orthonormal frame of $T_x \partial M$ at $x\in \partial M.$

The positive definite Laplace-Beltrami operator is 
\begin{equation*}
-\Delta_{g} u 
=  -\frac{1}{\sqrt{\det(g)}} \partial_{i} \left( g^{ij} \sqrt{\det(g)} \partial_{j} u \right).
\end{equation*}
For
\begin{equation*}
a = \frac{4(n - 1)}{n - 2}, p = \frac{2n}{n - 2}, 
\end{equation*}
the conformal Laplacian is 
\begin{equation*}
\Box_{g} u : = -a\Delta_{g} u + R_{g} u;
\end{equation*}
other sources define $\Box_g$ to be $a^{-1}$ times our $\Box_g.$
We use a  Robin boundary operator on $ \partial M $ as in \cite{E}:
\begin{equation}\label{Robin}
B_{g} u : = \frac{\partial u}{\partial \nu} + \frac{2}{p-2} h_{g} u.
\end{equation}
$ \eta_{1} = \eta_{1,g}$ is the first eigenvalue of conformal Laplacian with respect to the Robin condition $ B_{g} u = 0 $:
\begin{equation}\label{positive:eqn1}
\Box_{g} \varphi = \eta_{1} \varphi \; {\rm in} \; M, B_{g} \varphi = 0 \; {\rm on} \; \partial M.
\end{equation}
$\zeta_1$ is the first eigenvalue of the conformal Laplacian on the closed manifold $\partial M:$
$\Box_{i^*g}\psi = \zeta_1\psi$ on $\partial M.$

Let  $Y(\partial M, i^*g)$ be the Yamabe invariant on $\partial M$. We   have 
$$\Box_g>0 \Leftrightarrow \eta_1>0,\ \ \Box_{i^*g}>0 \Leftrightarrow \zeta_1>0 \Leftrightarrow Y(\partial M, i^*g)>0.$$
Here ${\rm sgn}(\zeta_1) = {\rm sgn}(Y (\partial M, i^*g))$ 
by \cite[Lem.~2.5]{KW}, and both signs are conformal invariants.

\subsection{The sign of 
$ R_{g} $ and $ R_{\imath^{*} g} $}

The prescribed scalar curvature problem both in the interior $ M $ and on the boundary $ \partial M $ is treated in \cite{XU1}. We recall that for given functions $ R \in \calC^{\infty}(M) $ and $ f \in \calC^{\infty}(\partial M) $, the existence of a metric $\tilde{g} \in [g]
$   with $R_{\tilde{g}} = R$ and   boundary condition $R_{\imath^*\tg} = f$ reduces to the existence of a positive, smooth solution of the Yamabe equation 
\begin{equation}\label{positive:eqn2}
\Box_{g} u = R u^{p-1} \; {\rm in} \; M, u = f \; {\rm on} \; \partial M.
\end{equation}
For an appropriate choice of $ f $,  
 the metric $\tilde g = u^{p-2} g = u^{4/(n-2)}g $ has  $ R_{\tilde{g}} = R, R_{\imath^{*} \tilde{g}} =f$.

Assuming the positivity of both $ \eta_{1} $ and $ \zeta_{1} $, we can improve the results in \cite{XU1} to control  the signs of both $ R_{\tg} $ and $ R_{\imath^{*} \tg} $.

\begin{theorem}\label{Scalar:thm1} 
Assume that $ \eta_{1} > 0 $ for $ (M^{n}, \partial M, g) $, $ n \geqslant 3 $. 
If either (i)  $ \zeta_{1} > 0 $ when $ n \geqslant 4 $, or 
(ii)  $ \chi(\partial M_i) > 0 $ for each boundary component $ \partial M_i $ of $\partial M$ when $ n = 3 $,
then there exists  $ \tilde{g} \in [g] $ such that $ R_{\tilde{g}} > 0 $ on $ M $ and $ R_{\imath^{*} \tilde{g}} > 0 $ on $ \partial M $.
\end{theorem}

\begin{proof}
Let $ \tilde{\eta}_{1} $ be the first eigenvalue of the conformal Laplacian $ \Box_{g} $ with vanishing Dirichlet boundary condition, {\it i.e.}, there exists $\varphi\in C^\infty(M)$ with
$ \Box_{g} \varphi = \tilde{\eta}_{1} \varphi \; {\rm in} \; M, \varphi = 0 \; {\rm on} \; \partial M.$   
By \cite[Prop.~4.3]{XU1}, $ \tilde{\eta}_{1} \geqslant \eta_{1} > 0 $. Given 
$ f \in \calC^{\infty}(\partial M) $, $f>0$,  there exists a positive, smooth function $ u $ and a positive constant $ \lambda $ such that 
\begin{equation}\label{Scalar:eqn2}
-a\Delta_{g} u + R_{g} u = \lambda u^{p-1} \; {\rm in} \; M, u = f \; {\rm on} \; \partial M,
\end{equation}
by \cite[Thm.~4.5]{XU1}.

By the solution of the Yamabe problem on closed manifolds,   $\zeta_1>0$ implies there exists $\phi\in C^\infty(\partial M), \phi >0$, 
such that the metric $g_1 =  \phi^{\frac{4}{n - 3}} \imath^{*} g $
has $R_{g_1} >0$ on $ \partial M $, provided  $ n \geqslant 4 $. When  $ n=3 $, there exists a conformal change 
$g_1 = e^{2\varphi} \imath^{*} g $ with Gaussian curvature $ R_{g_1 } > 0 $ on $ \partial M $
\cite{KW}.

Set
\begin{equation}\label{Scalar:eqn3}
f = \begin{cases} \left( \phi^{\frac{4}{n - 3}} \right)^{\frac{n - 2}{4}}, & n \geqslant 4, \\ \left( e^{2\varphi} \right)^{\frac{n - 2}{4}}, & 
 n = 3. \end{cases}
\end{equation}

Using this 
$ f $ in (\ref{Scalar:eqn2}), 
we conclude that $ \tilde{g} = u^{p-2} g $ has 
 $ R_{\tilde{g}} > 0 $ on $ M $ and $ R_{\imath^{*} \tilde{g}} > 0 $ on $ \partial M $.
\end{proof}
\medskip

\begin{corollary}\label{Ric:cor1}
Let $ (M^n, \partial M, g) $ have 
$\eta_1 >0$ and $ \zeta_{1} > 0$. Then there exists a  metric $ \tilde{g} \in [g] $ 
with $R_{\tilde g}$ a positive constant  
 on the interior of $ M $ and $R_{i^*\tilde g}$ a positive constant on $\partial M.$
\end{corollary}
\begin{proof} 
 Since $\eta_1>0,$ \cite[Thm.~1.1]{XU1} implies that 
there exists a metric $ \tilde{g} $ conformal to $ g $ with $R_{\tilde g}$ a positive constant  
 on the interior of $ M $ and $R_{i^*\tilde g}$ constant on $\partial M.$ On the closed manifold $\partial M$, it is well known that the sign of $\zeta_1$ equals the sign of the constant scalar curvature metric in $[g],$ so
$R_{i^*\tilde g}$ is a positive constant.   
\end{proof}

\begin{example}\label{ex:1} By \cite[Thm.~A]{Ba} and \cite[Thm.~4.5]{Hi}, every closed spin manifold $X$ of 
dimension $n\equiv 
0, 1, 3, 7\ ({\rm mod}\ 8)$ admits a metric $\bar g$ with nontrivial harmonic spinors. 
It is well known (see \S4.1)
that the spin bordism groups $\Omega^{\rm spin}_{2k-1}$ are torsion in odd dimensions, so outside the case
$n\equiv 0\ (\text{mod}\ 8)$, 
$NX$ ($N$ disjoint copies of $X$) is the boundary of a spin manifold $M$ for some 
$N$.
Then  $M$ admits no metric $g$ with $i^*g$ conformal to $\bar g$, such that $\eta_{1,g} >0$ and $g$ is a product metric near $NX$.
For in this case, the metric $\tg\in [g]$ with $R_{\tg}\equiv C$ has the constant $C$ positive, since 
$\sgn(C) = \sgn(\eta_{1,\tg}) = \sgn(\eta_{1,g}).$
For the product metric, 
we easily have $R_{\imath^*\tg} = R_{\tg}\equiv C>0,$ so there are no $\imath^*\tg$-harmonic spinors on $NX.$ Since the dimension of the space of harmonic spinors is a conformal invariant, there are no $\imath^* g$-harmonic spinors or $\bar g$-harmonic spinors, a contradiction.  This is an analytic version of the easier result that 
$M$ admits no psc metric which is a product near the boundary.
\end{example}

\subsection{Positivity of the conformal Laplacian:  boundary results}

We first note that Neumann boundary conditions are favorable in conformal geometry. 
 For $ \tilde{g} 
 = e^{2\phi} g $, the second fundamental form and mean curvature transform by
 \begin{align}\label{sff}
A_{\tilde{g}} (X, Y) &= e^{\phi} A_{g}(X, Y) + \frac{\partial \phi}{\partial \nu} e^{\phi} g(X, Y)\nu,\\
h_{\tilde{g}} &= e^{-\phi} h_{g} + \frac{\partial \phi}{\partial \nu} e^{-\phi}.\nonumber
\end{align}
Thus under this conformal change with Neumann boundary condition
$\frac{\partial \phi}{\partial \nu} = 0$, the vanishing of the second fundamental form or mean curvature for $g$ implies the same for $\tilde g$. Note that $ \frac{\partial \phi}{\partial \nu} = 0 $ is equivalent to $ \frac{\partial u}{\partial \nu} = 0 $ 
 for \begin{equation}\label{Revise0}
\tilde{g} = u^{p-2} g = e^{2\phi} g.
\end{equation}.

Assuming $\eta_1>0$ and good boundary geometry, we can alter a metric within its conformal class to have positive scalar curvature.

\begin{proposition}\label{APS:thm1}
 If $(M^{n}, \partial M, g) $   has $ \eta_{1} > 0 $ and $ \partial M $  totally geodesic, resp. minimal, 
  then there exists 
 a  metric $ \tilde{g}$ conformal to $g$  such that $ R_{\tilde{g}}>0 $ 
and $ \partial M $ remains  totally geodesic, resp. minimal, for $ \tilde{g} $.  
\end{proposition}
\begin{proof}
    We  treat the case where $\partial M$ is totally geodesic; the minimal boundary case is similar.
    By 
    \cite[\S1]{E}, the eigenvalue problem for the conformal Laplacian 
    \begin{equation}\label{Eigen1}
    -a\Delta_{g} \varphi  + R_{g} \varphi = \eta_{1} \varphi \; {\rm in} \; M, \frac{\partial \varphi}{\partial \nu} = 0 \; {\rm on} \; \partial M,
    \end{equation}
    admits a smooth solution $ \varphi $ which is positive on the interior of $M$. By the maximum principle, 
    $ \varphi > 0 $ on 
    $\partial M $ as well.  
    The mean curvature term $ h_{g} $ in the Robin boundary condition disappears,  since $ A_{g} \equiv 0 $.
    Rewriting (\ref{Eigen1}) as
    \begin{equation}\label{Eigen2}
    -a\Delta_{g} \varphi + R_{g} \varphi = (\eta_{1} \varphi^{2 - p} ) \varphi^{p-1} \; {\rm in} \; M, \frac{\partial \varphi}{\partial \nu}
    = 0 \; {\rm on} \; \partial M,
    \end{equation}
    it follows from (\ref{Eigen2}) that
    the metric $ \tilde{g} = \varphi^{p-2} g $ has scalar curvature $ R_{\tilde{g}} = \eta_{1} \varphi^{2 - p} > 0 $. 
        For $ \tilde{g} = \varphi^{p-2} g  
    := e^{\psi} g $, the boundary condition in (\ref{Eigen2}) gives
$\frac{\partial \psi}{\partial \nu} \equiv 0 \; {\rm on} \; \partial M.$
    Since $ A_{g} \equiv 0 $, (\ref{sff}) yields
$    A_{\tilde{g}}(X, Y)  = 0.$
 \end{proof} 
 
Similarly, we can reprove the known result that $ \eta_{1} > 0 $ implies the positivity of scalar curvature and vanishing of mean curvature for some $\tg\in [g].$
\begin{corollary}\label{APS:cor1}
 If $(M^{n}, \partial M, g) $ has $ \eta_{1} > 0 $, then there exists a  metric $ \tilde{g}\in [g]$  such that $ R_{\tilde{g}}>0 $ and $ h_{\tilde{g}} = 0 $.
\end{corollary}
\begin{proof}
    Again by  \cite[\S1]{E}, the eigenvalue problem for the conformal Laplacian 
    \begin{equation*}\label{Eigen3}
    -a\Delta_{g} \varphi  + R_{g} \varphi = \eta_{1} \varphi 
    = (\eta_{1} \varphi^{2 - p} ) \varphi^{p-1} \; {\rm in} \; M, \frac{\partial \varphi}{\partial \nu} + \frac{2}{p-2} h_{g} = 0 \; {\rm on} \; \partial M,
    \end{equation*}
    admits a smooth, positive solution $ \varphi $. As in (\ref{Eigen2}), 
   the metric $ \tilde{g} = \varphi^{p-2} g $ has the desired properties.
\end{proof}

\section{Positivity of the conformal Laplacian: A quasi-topological obstruction}

In this section, we find a ``quasi-topological" obstruction to a spin manifold with boundary admitting a metric with $\eta_1>0$ and $ \zeta_{1} > 0 $, namely an $\R$-valued $\xi$-invariant whose reduction mod $\Z$ is topological. As usual,
this is also an obstruction to psc metrics with minimal or constant positive mean curvature on the boundary.

Let $(M^n, g), n=4k$, be a  spin manifold with boundary $\partial M$, which is automatically spin.  
Let $S = S^+\oplus S^-$ be the spinor bundle on $M$, and let $\di= \di_{ g}:\Gamma(S^+)
\to \Gamma(S^-)$ be the Dirac operator on positive spinors on $M$. 
Let $(E=E_\alpha,\nabla)$ be a flat bundle with connection associated to a unitary representation $\alpha$ of $\pi_1(M)$.
The coupled Dirac operator  is
  $\di^E_{ g}=\di\otimes {\rm Id}+ {\rm Id}\otimes \nabla:\Gamma(S^+\otimes E)\to \Gamma(S^-\otimes E)$.

If $M$ has a metric $g$ with $\eta_1>0$,  
then we can conformally change $g$ to $
\tilde g$ so that $R_{\tilde g} >0 $ and $h_{\tilde g} = 0$ by Cor.~\ref{APS:cor1}.
 
\begin{lemma}\label{lem:4.1} Let  $(M,g)$ have  $\eta_1>0$, and let $\tilde g$ 
  be a conformally equivalent metric with $R_{\tilde g} >0 $ and $h_{\tilde g} = 0.$
  Then  $\ind(\di_{\tilde g}) =0$ for the Dirac operator with APS boundary conditions,
  and
    $\ind(\di^E_{\tilde g})=0 $ for the coupled Dirac operator with APS boundary conditions.
  \end{lemma}
  
  \begin{proof} The index $\ind(\di_{\tilde g})$ of the Dirac operator with APS boundary conditions vanishes by \cite[Thm.~4]{HMZ2}:   If $h_{\tilde g}\geqslant 0$, then the smallest eigenvalue $\lambda_1$ of $\di^E$ with APS boundary conditions satisfies 
$$\lambda_1^2 \geqslant \frac{n+1}{4n}\min_M R_{\tilde g}.$$
Since $R_{\tilde g} >0$, we cannot have $\lambda_1 = 0$.

  For the coupled Dirac operator, we just follow the proof of this theorem, replacing spinors $\psi\in \Gamma(S)$ by sections 
  $\psi\otimes s :=a^{ij}\psi_j\otimes s_i\in \Gamma(S\otimes E)$. 
  Specifically, we need equations (1) - (17) 
  and (the $\geqslant$ part of) Thm.~4 in 
  \cite{HMZ2} for $\psi\otimes s.$  
  We can assume $\{\psi_j\}, \{s_i\}$ are local orthonormal parallel sections of $S, E$, resp., at a fixed center point $x$ of a coordinate chart $U$,  with $a^{ij}\in C^\infty(U),$ by taking synchronous coordinates. Namely, we take an orthonormal frame $\{\psi_{j,x}\}$ of $S_x$ and parallel translate it out ({\it e.g.,} geodesic) rays starting at $x$ and filling out $U$.  Similarly, we parallel translate an orthonormal frame of $E_x$ along the rays to define
 $\{s_i\}$.

Let $\{e_k\}$ be an orthonormal frame of $T_xM.$ The key simplifications are  (dropping the $\tilde g$)
\begin{align*}\di^E(\psi\otimes s)&= 
\di(a^{ij}\psi_j)\otimes s_i
= e_k(a^{ij}) e_k\cdot \psi_j\otimes s_i,\\
\nabla_X^{S\otimes E}(\psi\otimes s) &:= (\nabla_X^S\otimes {\rm Id} + {\rm Id}\otimes \nabla_X^E)(
a^{ij}\psi_j\otimes s_i) = \nabla_X^S(a^{ij}\psi_j) \otimes s_i = 
X(a^{ij}) \cdot \psi_j\otimes s_i.  
\end{align*}
Here the dot is Clifford multiplication, and we use $\di \psi_j = e_k\nabla^S_{e_k}\psi_j = 0$ and $\nabla^Es_i=0$ at $x$.  For example,  for $X\in TM$,
\begin{align*} X\langle \psi\otimes s, \phi\otimes s'\rangle &=X\left(\langle a^{ij}\psi_j, b^{kl}
\psi_l\rangle_S \langle s_i,s_k\rangle_E\right) = X\langle a^{ij}\psi_j, b^{il}\psi_l\rangle_S\\
&= \sum_i \langle \nabla_X^S (a^{ij}\psi_j), b^{il}\psi_l\rangle +
\sum_i \langle  a^{ij}\psi_j, \nabla_X^S(b^{il}\psi_l)\rangle\\
&= \langle \nabla_X^{S\otimes E}(\psi\otimes s), \phi\otimes s'\rangle + 
\langle \psi\otimes s,  \nabla_X^{S\otimes E}(\phi\otimes s')\rangle,
\end{align*}
where the second line uses 
 \cite[(1)]{HMZ2}: $X\langle \psi, \phi\rangle_S = \langle \nabla^S_X \psi,\phi\rangle + \langle \psi, \nabla_X^S \phi\rangle.$  Thus we obtain the twisted version of  \cite[(1)]{HMZ2} from the untwisted version.  

Similarly, we get the twisted Lichnerowicz formula 
\begin{equation}\label{tl}(\di^E)^*\di^E = ( \nabla^{S\otimes E})^*
\nabla^{S\otimes E} + \frac{R}{4} \otimes {\rm Id} =\left( (\nabla^S)^*\nabla^S + 
\frac{R}{4}\right) \otimes {\rm Id} 
\end{equation}
 from the Lichnerowicz formula  $\di^*\di = (\nabla^S)^*\nabla^S + \frac{R}{4}.$  
 As a final example, in our notation \cite[(9)]{HMZ2} is
 \begin{equation}\label{9}\di^{S|_{\partial M}}\psi = \frac{n}{2}h\psi - \nu'\cdot \di^S\psi - \nabla^S_{\nu'}\psi,
 \end{equation}
 where $\nu'$ is the unit inward normal vector on $\partial M$.
 Then by (\ref{9}),
 \begin{align*}\di^{(S\otimes E)|_{\partial M}}(\psi\otimes s) &= 
 \di^{(S\otimes E)|_{\partial M}}(a^{ij}\psi_j\otimes s_i) = \di^{S|_{\partial M}}(a^{ij}\psi_j)\otimes s_i\\
 &=\left(\frac{n}{2}h a^{ij}\psi_j - \nu'\cdot \di^S(a^{ij}\psi_j) - \nabla^S_{\nu'}(a^{ij}\psi_j)\right)\otimes s_i\\
 &= \frac{n}{2}h (\psi\otimes s) -\nu'\cdot \di^{S\otimes E}(\psi\otimes s) -\nabla^{S\otimes E}_{\nu'} (\psi\otimes s).
 \end{align*}
 This is the twisted version of \cite[(9)]{HMZ2}.

 A straightforward but long check yields the twisted versions of (1)-(17), cumulating in the twisted version of \cite[Thm.~4]{HMZ2}: we still have
$\lambda_1^2 \geqslant [(n+1)/4n] \cdot\min_M R_{\tilde g}.$
 Thus the index of $\di^E_{\tilde g} $ vanishes.
    \end{proof}
  
We next show that 
the index of the Dirac operator with APS boundary conditions is a conformal invariant, whether or not $\eta_1>0.$
  \begin{lemma}\label{lem:4.1a} For $\tg$ conformal to $g$, $\ind(\di_g) = \ind(\di_{\tilde g}),$
  and $\ind(\di^E_g) = \ind(\di^E_{\tilde g}).$
\end{lemma}
 
\begin{proof}  
The APS theorem 
 for manifolds with non-product metric near the boundary  \cite[\S3]{Gilkey} is
\begin{equation}\label{APSG} \ind(\di_g) = \int_M \hat A(\Omega_g) +\int_{\partial M} P(g) -\frac{1}{2}({\rm dim\ ker}(\di_{i^*g}) 
+ \eta_g(0)),
\end{equation}
where $\hat A(\Omega_g)$ is the $\hat A$-polynomial in the curvature of $g$, and $P(g)$ is the Chern-Simons form related to the $\hat A$-polynomial.  As a polynomial in the Pontrjagin classes, $\hat A(\Omega_g)$ is a conformal invariant.  ${\rm Dim\ ker}(\di_{i^*g} )$ is a conformal invariant, since the Dirac operator is conformally covariant, {\it i.e.}. for $\tg = e^{2\phi} g$, we have $\di_{\tg} = e^{-(n+1)\phi/2} \di_g
e^{(n-1)\phi/2} $ \cite[\S4]{Hij1}, \cite[\S1.4]{Hi}.
The 
eta invariant is also a conformal invariant \cite{R}.  Finally, for $g_1$ a smoothing on $N = \partial M\times [0,1]$ of the metric $\imath^*g$ at $t=0$ to a product metric near $t=1,$
 we have $\int_{\partial M} P(g) = \int_N\hat A(\Omega_{g_1}).$  We can clearly construct $\tilde g_1$ conformal to $g_1$ on $N$ in the obvious notation.  As above, this implies that  $\int_{\partial M} P(g)$ is a conformal invariant. Thus $\ind(\di_g) = \ind(\di_{\tilde g})$.  
 
 For the twisted case, we check that
 \begin{align*}
\di_{\tg}^E(a^{ij}\psi_j\otimes s_i) &= (e^{-(n-1)\phi/2} \di_g
 e^{(n-1)\phi/2}\otimes {\rm Id} + {\rm Id} \otimes \nabla^E) (a^{ij}\psi_j\otimes s_i)\\
&= e^{-(n-1)\phi/2}e_k(e^{(n-1)\phi/2}a^{ij})e_k\cdot \psi_j\otimes s_i;
\end{align*}
\begin{align*}
\left(e^{-(n-1)\phi/2} \di_g^E e^{(n-1)\phi/2}\right)(a^{ij}\psi_j\otimes s_i)
&= e^{-(n-1)\phi/2} e_k(e^{(n-1)\phi/2}a^{ij})e_k\cdot\psi_j\otimes s_i.
 \end{align*}
 Thus the twisted Dirac operator is also conformally covariant. (This does not use that $\nabla^E$ is flat.)
 This implies that the $\eta$-invariant for $\di^E$ is also a conformal invariant.  As in the previous paragraph, $\ind(\di^E)$ is therefore a conformal invariant. 
 \end{proof}

Let $(Y,g)$ be a Riemannian spin manifold. The $\xi$-invariant associated to the unitary representation $\alpha$ of $\pi_1(Y)$  is
\begin{equation}\label{xidef}\xi_\alpha(Y, g) = \frac{1}{2}\left({\rm dim}\ \ker(\di^{E_\alpha}_{g}) + \eta^{E_\alpha}_g(0)\right) - \frac{k}{2}\left({\rm dim}\ \ker(\di_{g}) + \eta_g(0))
\right),
\end{equation}
where $\eta_g^{E_\alpha}(0)$  is the $\eta$-invariant for $\di_g^E$,
and $k$ is the rank of $\alpha.$ For a family of metrics $\{g_t\}$, $\xi_\alpha(Y, g_t)\in \R$ is not independent of $t$, but can only jump by an integer amount when the dimension of either the kernel of $\di_g$ or $\di^E_g$ jumps. Thus the mod $\Z$ reduction $\tilde \xi_\alpha(Y)$ of $\xi_\alpha(Y, g)$ 
 is a smooth invariant of $(Y,\alpha)$.  We call $\xi_\alpha(Y, g)$  quasi-topological, and now prove that it is an obstruction to $\eta_1>0, \zeta_1>0$ on spin manifolds.

\begin{theorem}\label{thm:4.1}  
Let 
$(M^{n},\partial M, g)$ be a spin manifold with boundary with $\eta_1>0$ and 
with $ \zeta_{1} > 0 $ on each component of $\partial M$. 
Assume that $\pi_1(\partial M)$ admits a unitary representation $\alpha$ that extends to a unitary representation of $\pi_1(M).$ Then $\xi_\alpha(\partial M, g) = 0.$ In particular, if $M$ admits a metric $\tg \in [g] $
with $R_{\tg}>0$, mean convex boundary, and $R_{\imath^*\tg}>0$, then  
$\xi_\alpha(\partial M, \imath^*g) = 0.$ 
\end{theorem}

\begin{proof}  
The left hand side of (\ref{APSG}) vanishes for  $\tg$ as in
Lem.~\ref{lem:4.1}, and hence for $g$ by Lem.~\ref{lem:4.1a}. 
  On the right hand side, $P( g) = P(\Omega_{ g}, \omega_{ g})$ is a
an invariant polynomial in the components of the curvature and connection forms of $g$ \cite[\S3]{Gilkey}.  
By  
Thm.~\ref{Scalar:thm1} and the Lichnerowicz formula, there 
exists $\tg\in [g]$ with
ker$(\di_{i^*\tg})=0$ on $\partial M$.  By the conformal invariance of this kernel, 
ker$(\di_{i^* g})=0$. 
Therefore,
\begin{equation}\label{APS2}\eta_{ g}(0) = -2\int_M \hat A(\Omega_{ g}) -2\int_{\partial M} P( g)
\end{equation}
is the integral of  locally computed curvature and connection terms.

The APS formula for $\di_g^E$ is
\begin{align*}\ind(\di^E_g) &= \int_M \hat A(\Omega_g)ch(\Omega_E) +\int_{\partial M} P^E(g) -\frac{1}{2}({\rm dim\ ker}(\di^E_{i^*g}) 
+ \eta^E_g(0))\\
&= k\int_M \hat A(\Omega_g) + k\int_{\partial M}P(g) -\frac{1}{2}({\rm dim\ ker}(\di^E_{i^*g}) 
+ \eta^E_g(0)).
\end{align*}
This follows from  
\begin{align*}P^E(g) &= c \int_0^1 \Tr (\omega^{S\otimes E}_t, \Omega_t^{S\otimes E},\ldots, \Omega_t^{S\otimes E})dt
= c \int_0^1 \Tr (\omega^S_t\otimes {\rm Id},\Omega^S_t\otimes {\rm Id},\ldots,\Omega^S_t\otimes {\rm Id})dt\\
&= ck\int_0^1 \Tr (\omega^S_t,\Omega^S_t,\ldots,\Omega^S_t)dt
\end{align*}
where $c$ is a dimension constant, in the usual notation for Chern-Simons forms.  As in the last paragraph, we can assume that $\ind(\di^E_g)=0$.  The argument in 
Thm.~\ref{Scalar:thm1}  applies to $\di^E_{i^*g}$, by the twisted Licherowicz formula (\ref{tl}).  Thus we can also assume
${\rm dim\ ker}(\di^E_{i^*g})=0$.  We obtain
\begin{equation}\label{APS3}\eta^E_{ g}(0) = -2k\int_M \hat A(\Omega_{ g}) -2k\int_{\partial M} P( g).
\end{equation}
From (\ref{APS2}), (\ref{APS3}), we get 
\begin{align*}\xi_\alpha(\partial M, g)& = 
\frac{1}{2} \eta^{E_\alpha}_g(0) - \frac{k}{2} \eta_g(0)
\\
&= -k\left(\int_M \hat A(\Omega_{ g}) +\int_{\partial M} P( g)\right) +k
\left(\int_M \hat A(\Omega_{ g}) +\int_{\partial M} P( g)\right)\\
&=0.
\end{align*}

The last statement follows from the variational definitions of $\eta_{1,\tg}$ and $\zeta_{1,\imath^*\tg}.$
\end{proof}

\begin{remark} (i) This Theorem implies that the nonvanishing of the topological invariant $\tilde\xi_\alpha(Y)$ is an obstruction to $\eta_1>0$, but in practice computing $\tilde\xi_\alpha(Y)$ comes down to computing $\xi_\alpha(Y,g)$ and reducing mod $\Z.$

(ii) (\ref{APS2}) is proved in \cite[Thm.~3.9]{APSII},  if $M$ admits a psc metric $ g$ which is a product near $\partial M.$  Such a metric has vanishing second fundamental form, so $h_{ g} = 0$.  Thus these metrics have 
$\eta_1>0$, and so form a special case of the Theorem.
\end{remark}

\section{Manifolds with \texorpdfstring{$3$}--dimensional lens space boundaries}
There are explicit formulas for the $\xitilde$-invariants of spherical space forms; these are particularly tractable for lens spaces with odd order fundamental group. Making use of these formulas, we present some concrete examples of spin $4$-manifolds with boundary a disjoint union of lens spaces that cannot admit a metric with $\eta_1>0$ and $\zeta_1>0$; geometric restrictions are given in Thm.~\ref{thm:5.1}.

The lens space 
$L(p,q)$ is the quotient of 
$S^3\subset \C^2$ by the action of $\Z_p$ generated by 
\[ 
T \cdot(z_1,z_2) = (\omega z_1,\omega^q z_2),\ \text{where}\ \omega = e^{2\pi i/p}.\]
Here $T$ is a fixed generator of $\Z_p$, which can also be viewed as a fixed generator of $\pi_1(L)$.  The round metric on $S^{3}$ descends to a metric on $L$, which evidently has positive scalar curvature. 
The $\xitilde$-invariants with respect to this standard metric have been computed by Gilkey. The precise formula is a little unwieldy, but we only need a non-vanishing result. 
\begin{lemma}[Gilkey]\label{L:lens}
Let $\alpha: \pi_1(L(p,q)) \to \U(1)$ take $T$ to $\omega$. Then, with respect to any metric $g$ of positive scalar curvature, $\xi_\alpha(L(p,q),g) \neq 0$.
\end{lemma}
\begin{proof}
For $p$ odd (so there is a unique spin structure), Gilkey~\cite[Thm.~2.5(c)]{gilkey:eta-odd} shows that with respect to the standard round metric $g_0$ on $L(p,q)$
    \[
\xi_{\alpha}(L(p,q), g_0)\; =\; - (d/p) \cdot (p + 1)/2 \pmod \Z,
\]
where $d$ is a certain integer relatively prime to $48$ and $p$. One can easily see that the right hand side is never zero modulo the integers, and hence  $\xi_{\alpha}(L(p,q),g)$ is not zero as a real number.

A recent result of Bamler and Kleiner~\cite{bamler-kleiner:contractible} shows that the space of psc metrics on any $3$-manifold is contractible, and in particular is path-connected. By the vanishing of the kernel of the Dirac operator on a manifold with a psc metric, there is no spectral flow for the $\eta$-invariant of the Dirac operator along a path of psc metrics. The same holds for the twisted Dirac operator, and so $\xi_\alpha$ is constant along a path of psc metrics. It follows that $\xi_{\alpha}(L(p,q),g)$ is non-zero for any psc metric $g$ on $L(p,q)$.
\end{proof}

\subsection{A rational coboundary for lens spaces}\label{s:lens}
An argument with the Atiyah-Hirzebruch spectral sequence (AHSS)~\cite{davis-kirk:at} says that the bordism group $\Omega^{\spin}_{2k+1}(BG)$ is torsion for any finite group $G$. In other words, for any odd-dimensional spin manifold $Y$ and unitary representation $\alpha$ of $\pi_1(Y)$ that factors through a finite group, there is a spin manifold $M$ with a unitary representation $\hat{\alpha}$ of its fundamental group, and an integer $N$ with $\partial(M,\hat{\alpha}) = N (Y,\alpha)$. 

Here is a brief outline of the argument. The AHSS has $E^2_{p,q} \cong H_p(BG;\Omega^{\spin}_{q})$ and converges to $\Omega^{\spin}_{*}(BG)$. It is standard that for a finite group, the rational homology of $BG$ vanishes in degree greater than zero.  Moreover, the Anderson-Brown-Peterson calculation of the spin cobordism group~\cite{abp:spin} (see~\cite[p.~336]{stong:cobordism} for the rational calculation) implies that $\Omega^{\spin}_{q} \otimes \Q = 0$ when $q$ is odd. It follows that the AHSS for $\Omega^{\spin}_{*} \otimes \Q$  collapses immediately at the $E^2$ term, and so $\Omega^{\spin}_{2k+1}(BG)\otimes \Q=0$.

For the lens spaces with odd order fundamental group, we can give an explicit construction. The multiple $N$ 
we find is not optimal; for instance, 
we find a spin $4$-manifold $M$ with $\partial(M,\hat\alpha) = p^2(L(p,q),\alpha)$, whereas a closer examination of the AHSS shows that $\Omega^{\spin}_{3}(B\Z_p) \cong \Z_p$ so $N = p$ would do. 

\begin{figure}[ht]
\labellist
\small\hair 2pt
\pinlabel{$\tau$} [ ] at 115 387
\endlabellist
\centering
\includegraphics[scale = .4]{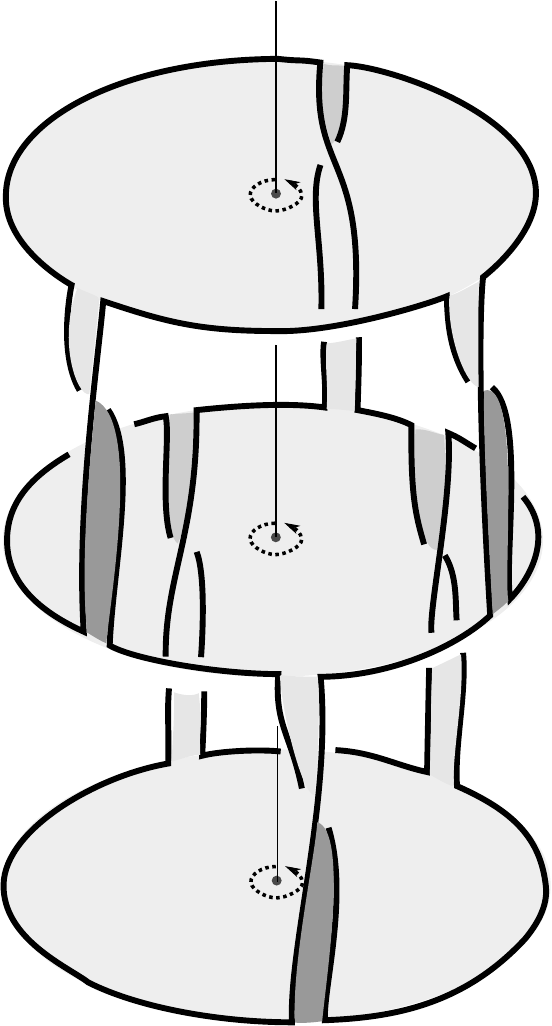}
  \caption{The surface $\Sigma_3^0$}\label{F:Sigma}
\end{figure}

The construction begins with the surface $\Sigma_p^0$ drawn (for $p=3$) in Figure~\ref{F:Sigma}, inspired by 
Figure 1 in~\cite{gordon:G-signature} and~\cite[\S 5]{gilmer:thesis}. 
It may described as a stack of $p$ discs, joined by $p$ half-twisted bands as shown in the figure. The surface $\Sigma_p^0$ has a $p$-fold symmetry $\tau$, given by counterclockwise rotation by angle $2\pi/p$ with fixed point set $p$ points, one on each of the discs. Furthermore $\partial \Sigma_p^0$ is $p$ circles, cyclically permuted by $\tau$. Hence the symmetry extends to an action of $\Z_p$ on the closed surface 
\[
\Sigma_p = \Sigma_p^0 \cup p D^2.
\]

Note that the symmetry $\tau$ on $\Sigma_p^0$ comes from the $p$-fold symmetry of $S^3$ given by rotation around the indicated axis. Now $S^3$ is a spin manifold, and it is standard that this symmetry lifts to the spin bundle of $S^3$. Therefore $\tau$ lifts to the spin bundle of the restriction of the spin structure on $S^3$ to $\Sigma_p^0$. 

We claim that this equivariant spin structure $\spincs$ extends to $\Sigma_p$ as well. 
To see this, it suffices to show that the induced spin structure on each boundary component of $\Sigma_p^0$ is the bounding spin structure on $S^1$ (and hence extends, uniquely, over each of the $p$ discs that are added.). By construction, $p$ copies of the spin structure on the circle is a boundary, for it bounds $\Sigma_p^0$. But the spin cobordism group $\Omega_2^{\spin}$ is $\Z_2$, and since $p$ is odd, it follows that $(S^1,\spincs) = 0$ in $\Omega_2^{\spin}$. Alternatively, one can note that the normal framing on each boundary component induced by the inward normal of $\Sigma_p^0$ is $p-1$, which is even. Hence the $4$-manifold obtained by adding $2$-handles to $S^3 \times I$ along $\partial \Sigma_p^0$ with framing $p-1$ is spin, and contains $\Sigma_p$. The spin structure on this $4$-manifold induces a spin structure on $\Sigma_p$ extending the given one on $\Sigma_p^0$.

To construct a bounding manifold  for  
 $p^2L(p,q)$, 
consider the product $(\Sigma_p)^2$ of two copies of $\Sigma_p$. There is an action of $\Z_p$ on $(\Sigma_p)^2$, where the generator, say $T$, acts by $\tau$ on the first factor of $\Sigma_p$ and by $\tau^q$ on the second factor. The fixed point set of $T$ therefore consists of $p^2$ isolated points, where the local representation is exactly the action that defines the lens space $L(p,q)$. Now remove an invariant $D^{4}$ neighborhood of each fixed point, and take the quotient by $\Z_p$. The quotient is a spin manifold with boundary $p^2 L(p,q)$. The orientation on the lens space boundaries is opposite to the orientation implicit in the definition, so we  set $M= M(p,q)$ to be this spin manifold with orientation reversed.  By construction, there is a homomorphism $\hat\alpha: \pi_1(M) \to \Z_p$ that restricts to $\alpha$ on each copy of $L(p,q)$. 

Combining Thm.~\ref{thm:4.1} with the construction above gives our main example. Recall the definition of
$\tg\in [g]$, $k\in [psc], k'\in [psc,mc] $ in the Introduction.
\begin{theorem}\label{thm:5.1} 
There is no metric $g$ on $M(p,q)$ with $\eta_{1,g} >0$ and $\zeta_{1,\imath^*g}>0. $ 
 In particular, (i) for $k\in [psc]$ on $p^2L(p,q)$, 
  there is no metric $g$ on  $M(p,q)$  with $\eta_{1,g} >0$ and $\imath^*g\in [k]$; 
  (ii)
no metric $k\in [psc]$ on $p^2L(p,q)$  
extends to a  metric $k'\in [psc,mc]$ on $M(p,q)$;
(iii) statement (ii) holds if we replace ``mean convex boundary" with ``minimal boundary" or ``is a product near the boundary."  
\end{theorem}

\begin{proof} By Thm.~\ref{thm:4.1}, $\eta_{1,g} >0,\zeta_{1,\imath^*g}>0$ 
implies $\xi_{\alpha}(p^2L(p,q), \imath^{*} g) = 0. $
This contradicts Lem.~\ref{L:lens}.

(i) follows from the main statement, since a  metric $k\in [psc]$ has $\zeta_{1,k}>0.$ Indeed,
the conformal covariance of $\Box_k$ implies that the sign of $\zeta_1$ is a conformal invariant.
For (ii), we know $\zeta_{1,k} >0$. Since the sign of the relative Yamabe invariant is a conformal invariant, Lem.~\ref{prop:RYI} below implies the sign of $\eta_{1,k'}$ is a conformal invariant. A psc metric $g$ with mean convex boundary has $\eta_{1,g} >0.$  Thus $\eta_{1,k'}>0,$ so we get a contradiction as above.
 Finally,  (ii) implies (iii).
 \end{proof}

 \begin{remark}  
 If we 
replace the hypothesis with $\eta_1>0$ and minimal boundary,
the proof follows from earlier work   \cite[Thm.~2.2]{AB2}, \cite[Thm.~4.6]{AB1}, where it is proved that such a metric can be conformally deformed to a product metric near the boundary.
Since psc on $M(p,q)$ plus minimal boundary implies $\eta_{1}>0,$ our results and the results of Akutagawa-Botvinnik generalize the case of a psc metric which is a product near the boundary, the case originally treated in
\cite{APSI}.  However, in \S5 we prove that all these results can be reduced to this product metric case.
 \end{remark}

\subsection{Higher dimensions}  

Let $(Y,g,\alpha)$ be a  $(4k-1)$-dimensional Riemannian spin manifold with $k>1$, where $\alpha$ is a representation of $\pi_1(Y)$ that factors through a finite group.
We know $N(Y, \alpha) = \partial (M, \hat\alpha)$ for some $N$ with $M$ spin 
and $\hat\alpha$ an extension of $\alpha$. 
If $\xi_\alpha(Y,\imath^*g')\neq 0$ for some metric $g'$ on $M$, then  $M$ admits no metric $g$ with $\eta_{1,g}>0$ and $\zeta_{1,i^*g}>0$ on each component of $\partial M.$  As above, the $\xi$-invariant is an obstruction to extending $k\in [psc]$ from $NY$ to a metric $k'\in [psc, mc]$ on a bounding manifold $M$. 

In the quasi-topological direction,
we list two cases where $M$ admits no metric with $\eta_1 >0, \zeta_1>0$, but we know of no examples.

\begin{enumerate}
\item Assume  $\xitilde_\alpha(Y) \neq 0\in \R/\Z,$ so  $\xi_\alpha(Y,g) \neq 0\in \R.$    As in Cor.~\ref{Ric:cor1}, we may assume that $i^*g$ is psc.   If the space of psc metrics is connected, then as in Lem.~\ref{L:lens}, $\xi_\alpha(NY,g)=N\xi_\alpha(Y,g)= 0.$  This contradicts $\eta_1>0, \zeta_1>0$, by Thm.~\ref{thm:4.1} .  There are examples of higher dimensional lens spaces with nonzero $\xi$-invariants, but we don't have an example of a higher dimensional manifold with connected space of psc metrics.

\item If the space of psc metrics is not connected, then for $NY = \amalg_{i=1}^N Y_i$,\\
$\xi_\alpha(NY,g) = \sum_{i=1}^N \xi_\alpha(Y_i,g)$ may vanish, since the $\xi_\alpha(Y_i,g)$ differ by integers.  However, this cannot occur if $\xi_\alpha(Y_i,g)\not\in \mathbb Q.$  While irrational $\xi$-invariants occur for $Y= S^1$, we again know of no higher dimensional examples.  In particular, the examples in \cite[Lem.~2.3]{BG} are all rational.

\end{enumerate}

\section{Positivity of the conformal Laplacian and the relative Yamabe invariant}

In this section, we relate the positivity of the conformal Laplacian to the positivity of the relative Yamabe invariant in \cite{AB1}.  Using  powerful techniques from \cite{AB1}, \cite{LZ}, this allows us to reprove the main Theorems~\ref{thm:4.1}, \ref{thm:5.1} more quickly, 
by reducing these theorems to the case where the psc metric is a product metric near the boundary.  While these proofs are shorter, they are less elementary and lose some of the results proved in \S\S2,3.

To define the relative Yamabe invariant, let $M$ be an oriented manifold with boundary $\partial M$, and let $k$ be a Riemannian metric on $\partial M$.
Let $\bar\calC$ be a conformal class of metrics on $M$, and set $\partial \bar\calC = [k]$ if  $\bar g\in \bar\calC$ has $\imath^* \bar g \in [k]$, for the inclusion $i:\partial M\to M.$  Set
$\bar\calC^o =\{\bar g\in \bar\calC: h_{\bg} = 0\}.$
Each $\bar\calC$ contains a metric $\bg_0 
\in \bar\calC^0$, 
so $$\bar\calC^o =\{ u^{4/(n-2)}\bg_{0}
: u>0, u\in C^\infty, \partial_\nu u = 0\}.$$
By definition, the relative Yamabe invariant is
\begin{align}\label{eq:RYI}
Y(M,L,[k]) &:= \sup_{\bar\calC,\partial \bar \calC = [k]} Y_{\bar\calC}(M,L,[k])
:= \sup_{\bar\calC,\partial \bar \calC = [k]}\inf_{\bg\in \bar\calC^o}\frac{\int_M R_{\bg}\ d\text{Vol}_{\bg}}{\text{ Vol}_{\bg}( M)^{\frac{n-2}{n}}}.
\end{align}

The following Proposition is well-known for closed manifolds, but needs a modified proof for manifolds with boundary.
\begin{proposition}\label{prop:RYI}
$$\eta_{1,g}>0 \Leftrightarrow Y(M,\partial M,[\imath^*g])>0.$$
\end{proposition}

\begin{proof}
Assume that $\eta_{1,g}>0.$
For  $\bg = u^{4/(n-2)}g \in \bar\calC$ 
the ratio in the last term in (\ref{eq:RYI}) is the same (up to a positive constant $A$) as 
$$\frac{\int_{M} \left(\frac{4(n - 1)}{n - 2} \lvert \nabla_{g} u \rvert^{2}  +  R_{g} u^{2}\right) \dvol
+ 2(n - 1)(n - 2)\int_{\partial M} h_g u^2d\text{Vol}_{\imath^*g}}{\left( \int_{ M} u^{\frac{2n}{n - 2}} \dvol \right)^{\frac{n - 2}{n}}},$$
by the conformal transformation formulas for scalar curvature and volume.  
Thus $Y(M,\partial M,[k])>0$ if there exists $\bar\calC$ with $\partial\bar\calC = [k]$ such that 
\begin{align} \label{eq:111}
0 &< 
\inf_{ \{  u\in C^\infty: u>0, \partial_\nu u = 0\}}\frac{\int_{M} \left(\frac{4(n - 1)}{n - 2} \lvert \nabla_{g} u \rvert^{2}  +  R_{g} u^{2}\right) \dvol}
{\left( \int_{M} u^{\frac{2n}{n - 2}} \dvol \right)^{\frac{n - 2}{n}}} := \inf_{ \{  u\in C^\infty: u>0, \partial_\nu u = 0\}}Q(u),
\end{align}
for $g\in \bar\calC$ with $h_g=0.$  Note that the Robin boundary condition $B_{g} u  = \frac{\partial u}{\partial \nu} + \frac{2}{p-2} h_{g} u=0$ is satisfied for $u\in \bar\calC^o.$

By Cor.~\ref{APS:cor1}, we may assume $\eta_{1,g} >0$ for some $g\in \bar\calC$ with $R_g>0, h_g=0$. Then
\begin{align}\label{eq:112}0&< \eta_{1,g} = \inf_{\{u\in C^\infty: \partial_\nu u + h_g u=0\}}\frac{\int_{M} \left(\frac{4(n - 1)}{n - 2} \lvert \nabla_{g} u \rvert^{2}  +  R_{g} u^{2}\right) \dvol
}{ \int_{M} u^2 \ \dvol }\\
&< \inf_{\{u\in C^\infty: u>0, \partial_\nu u =0\}}\frac{\int_{M} \left(\frac{4(n - 1)}{n - 2} \lvert \nabla_{g} u \rvert^{2}  +  R_{g} u^{2}\right) \dvol
}{ \int_{M} u^2 \ \dvol }.\notag
\end{align}
To compare the denominators in (\ref{eq:111}) and (\ref{eq:112}), we use \cite[Thm.~0.2]{LZ}:
there exists a constant $ A_{1} = A_{1}(M, g) > 0 $ such that for all $ u \in H^{1}(M, g) $,
\begin{equation*}
    \left(\int_{M} u^{\frac{2n}{n - 2}} \dvol \right)^{\frac{n - 2}{n}} \leqslant 2^{\frac{2}{n}} S_{1} \int_{M} \lvert \nabla_{g} u \rvert^{2} \dvol + A_{1} \int_{M} u^{2} \dvol + A_{1} \int_{\partial M} u^{2}
    d\text{Vol}_{\imath^*g}.
\end{equation*}
($ S_{1} $ is the Sobolev constant for the upper half plane $ \mathbb{R}_{+}^{n} $.)  Therefore,
\begin{align*}
Q(u) 
& \geqslant \frac{\int_{M} \left(\frac{4(n - 1)}{n - 2} \lvert \nabla_{g} u \rvert^{2}  +  R_{g} u^{2} \right) \dvol
}{ 2^{\frac{2}{n}} S_{1} \int_{M} \lvert \nabla_{g} u \rvert^{2} \dvol + A_{1} \int_{M} u^{2} \dvol + A_{1} \int_{\partial M} u^{2} d\text{Vol}_{\imath^*g} } \\
& \geqslant \frac{\int_{M} \left(\frac{4(n - 1)}{n - 2} \lvert \nabla_{g} u \rvert^{2}  +  R_{g} u^{2} \right) \dvol
}{ 2^{\frac{2}{n}} S_{1} \int_{M} \lvert \nabla_{g} u \rvert^{2} \dvol + A_{1} \int_{M} u^{2} \dvol + A_{1}C ( \int_{M} \lvert \nabla_{g} u \rvert^{2}\dvol + \int_{M} u^{2} \dvol)},
\end{align*}
where we use the trace theorem on the last term in the denominator:  
$$ \lVert u \rVert_{L^{2}(\partial M)} \leqslant C \lVert u \rVert_{H^{1}(M)} = C\left(\int_{M} \lvert \nabla_{g} u \rvert^{2}\dvol + \int_{M} u^{2} \dvol\right). $$
By the positivity of $R_g$, there exist positive constants $B_1, B_2$ such that 
\begin{equation*}
    Q(u) \geqslant \frac{B_{2} \left( \int_{M} \lvert \nabla_{g} u \rvert^{2} \dvol + \int_{M} u^{2} \dvol \right)}{B_{1} \left( \int_{M} \lvert \nabla_{g} u \rvert^{2} \dvol + \int_{M} u^{2} \dvol \right)} = \frac{B_{2}}{B_{1}}, \forall u \in H^{1}(M).
\end{equation*}
By (\ref{eq:111}), $Y(M,\partial M,[k]) >0.$
\medskip

Now assume $Y(M,\partial M,[\imath^*g])>0.$ By H\"older's inequality, there is a positive constant $B$ such that 
\begin{equation*}
    \int_{M} u^2 \ \dvol \leqslant B \left( \int_{M} u^{\frac{2n}{n - 2}} \dvol \right)^{\frac{n - 2}{n}}
\end{equation*}
As above, this implies
\begin{equation*}
 0<  \inf_{u\in C^\infty}  \frac{\int_{M} \left(\frac{4(n - 1)}{n - 2} \lvert \nabla_{g} u \rvert^{2}  +  R_{g} u^{2} \right) \dvol
}{ (\int_{M} u^\frac{2n}{n - 2} \ \dvol)^{\frac{n-2}{n}} } \leqslant B\inf_{u\in C^\infty}\frac{\int_{M} \left(\frac{4(n - 1)}{n - 2} \lvert \nabla_{g} u \rvert^{2}  +  R_{g} u^{2} \right) \dvol
}{ \int_{M} u^2 \ \dvol }.
\end{equation*}
Therefore, $\eta_{1,g}>0.$
\end{proof}

We can now relate our spectral positivity hypothesis to geometric conditions.  

\begin{corollary}\label{cor:1}  
There exists a metric $g$ on $M$ with $\eta_{1,g}>0, \zeta_{1,\imath^*g}>0$ on $M$ and $\partial M$ iff there exists a psc metric $k$ on $\partial M$ which extends to a psc metric $k'$ on $M$ which is a product near $\partial M$.
\end{corollary}

\begin{proof}
$\Rightarrow)$   Since $\zeta_{1,\imath^*g}>0$, using the associated eigenfunction $\phi$ as in the proof of Cor.~\ref{APS:cor1}  yields a
a psc metric $k\in [\imath^*g]$.
 Because $\eta_{1,g}>0$, by Prop.~\ref{prop:RYI}, $Y(M,\partial M,[\imath^*g] = [k])>0.$ By \cite[Cor. B]{AB1}, $k$ extends to a psc metric $k'$ on $M$ which is a product near the boundary.  

$\Leftarrow)$  Because $k$ is psc, $\zeta_{1,k} >0.$ Because $k'$ exists, $Y(M, \partial M, [k])>0$ by 
\cite[Cor. B]{AB1}. By Prop.~\ref{prop:RYI}, $\eta_{1,k'} >0.$
\end{proof}

This corollary leads to quick proofs of the main statements in the main theorems.
\medskip

\noindent {\bf Theorems 3.1 and 4.1}
(i) Let 
$(M^{n},\partial M, g)$ be a spin manifold with boundary with $\eta_1>0$ and 
with $ \zeta_{1} > 0 $ on each component of $\partial M$. 
Assume that $\pi_1(\partial M)$ admits a unitary representation $\alpha$ that extends to a unitary representation of $\pi_1(M).$ Then $\xi_\alpha(\partial M, g) = 0.$

(ii) There is no metric $g$ on $M(p,q)$ with $\eta_{1,g} >0$ and $\zeta_{1,\imath^*g}>0. $ 
\bigskip

\begin{proof}  (i) By Cor~\ref{cor:1} we may assume that $g$ is a psc metric on $M$ with product structure near 
$\partial M$.  Then by \cite[Thm.~3.9]{APSI} and the following discussion, we get 
$\xi_\alpha(\partial M, g) = 0.$

(ii) As \S4, this follows immediately from (i) and Lem.~\ref{L:lens}
\end{proof}

\bibliographystyle{plain}
\bibliography{DannySteveFrank}
\end{document}